\documentclass[11pt]{amsart}

\usepackage[utf8]{inputenc}
\usepackage[T1]{fontenc}
\usepackage{amsmath,amssymb,amsthm,mathtools,mathrsfs}
\usepackage[margin=1.1in]{geometry}
\usepackage[dvipsnames]{xcolor}
\definecolor{srcA}{RGB}{0,84,166}
\definecolor{srcX}{RGB}{150,95,0}
\definecolor{srcC}{RGB}{120,40,140}

\definecolor{bodylinkblue}{RGB}{0,64,255}
\definecolor{toclinkred}{RGB}{255,0,0}
\definecolor{tocsubsectionred}{RGB}{255,0,0}
\usepackage[colorlinks=true,linkcolor=toclinkred,citecolor=bodylinkblue,urlcolor=bodylinkblue]{hyperref}

\numberwithin{equation}{section}
\allowdisplaybreaks
\makeatletter
\def\l@section{\@tocline{1}{0pt}{1pc}{}{\color{toclinkred}}}
\def\l@subsection{\@tocline{2}{0pt}{1pc}{5pc}{\color{tocsubsectionred}}}
\def\l@subsubsection{\@tocline{3}{0pt}{1pc}{7pc}{\color{tocsubsectionred}}}
\renewcommand{\tocsection}[3]{%
  \textcolor{toclinkred}{\indentlabel{\@ifnotempty{#2}{\ignorespaces#1 #2.\quad}}#3}}
\renewcommand{\tocsubsection}[3]{%
  \textcolor{tocsubsectionred}{\indentlabel{\@ifnotempty{#2}{\ignorespaces#1 #2.\quad}}#3}}
\renewcommand{\tocsubsubsection}[3]{%
  \textcolor{tocsubsectionred}{\indentlabel{\@ifnotempty{#2}{\ignorespaces#1 #2.\quad}}#3}}
\def\@settitle{\begin{center}%
  \baselineskip14\p@\relax
  \bfseries
  \uppercasenonmath\@title
  \@title
  \end{center}%
}
\def\@setauthors{%
  \begingroup
  \def\thanks{\protect\thanks@warning}%
  \trivlist
  \centering\footnotesize \@topsep30\p@\relax
  \advance\@topsep by -\baselineskip
  \item\relax
  \author@andify\authors
  \def\\{\protect\linebreak}%
  \MakeUppercase{\authors}%
  \ifx\@empty\contribs\else
    ,\penalty-3 \space \@setcontribs
    \@closetoccontribs
  \fi
  \endtrivlist
  \endgroup
}
\def\@setabstracta{%
  \ifvoid\abstractbox\else
    \skip@38\p@ \advance\skip@-\lastskip
    \advance\skip@-\baselineskip \vskip\skip@
    \box\abstractbox
    \prevdepth\z@
  \fi
}
\makeatother

\newcommand{\D}{\mathbb D}
\newcommand{\w}{\omega}
\newcommand{\C}{\mathbb C}

\newcommand{\abs}[1]{\lvert #1\rvert}

\newcommand{\mean}[1]{\langle #1\rangle}
\DeclareMathOperator{\Ric}{Ric}

\DeclareMathOperator{\Vol}{Vol}

\DeclareMathOperator{\Hess}{Hess}
\newcommand{\dd}{\,d}
\newcommand{\eps}{\varepsilon}
\newcommand{\Rz}[1]{R_\w(#1)}

\newtheorem*{brody}{Brody's criterion}
\newtheorem{mainthm}{Theorem}

\newtheorem{theorem}{Theorem}[section]
\newtheorem{lemma}[theorem]{Lemma}
\newtheorem{proposition}[theorem]{Proposition}

\newtheorem{question}[theorem]{Question}

\theoremstyle{definition}

\newtheorem{example}[theorem]{Example}
\theoremstyle{remark}
\newtheorem{remark}[theorem]{Remark}

\title[Remark on Kobayashi's conjecture: a volume Schwarz lemma]{Remark on Kobayashi's conjecture: a volume Schwarz lemma}

\author{Shiyu Zhang}

\date{19 November 2026}

\begin{document}
\enlargethispage{4pt}
\begin{abstract}
This remark offers an expository reconstruction of OpenAI’s proof of Kobayashi’s ampleness conjecture for compact K\"ahler manifolds, together with new results. We extract a volume Schwarz lemma for extremal discs and derive a pointwise volume comparison: the volume form of a K\"ahler metric with $\mathrm{Ric}(\omega)\ge-\omega$ is bounded below by an intrinsic volume form built from the K\"ahler--Einstein metrics of pseudoconvex domains containing the Kobayashi--Royden indicatrix. As a consequence, if a compact K\"ahler manifold $X$ has nef canonical bundle and its Kobayashi--Royden indicatrix is bounded at some point, then $X$ is of general type.
\end{abstract}
\maketitle

\section{Introduction}

Throughout this remark, hyperbolic means Kobayashi hyperbolic.

\subsection{Background}\label{sec:background}
A complex manifold is hyperbolic if its Kobayashi pseudo-distance is a distance. Kobayashi asked the following question \cite[Problem~6]{Kob}, \cite[Problem~F.1]{Kob76}.
\begin{question}
 Is a compact hyperbolic complex manifold projective and of general type?
\end{question}
Lang made a stronger prediction \cite[Conjecture~7.2]{Lang86} that the canonical bundle of a projective hyperbolic manifold is ample. This is now known as Kobayashi's conjecture \cite{Pet91,Cam91}. There has been significant progress in some special cases:
it holds for compact manifolds satisfying Gromov's K\"ahler hyperbolicity
\cite{CY}, and Wu-Yau's theorem for compact K\"ahler manifolds with
negative holomorphic sectional curvature \cite{WY,TY} (see \cite{WY16remark,DT} for some extensions and \cite{DiverioSurvey} for a survey). Assuming abundance, the
conjecture reduces to excluding hyperbolic Calabi--Yau manifolds
(see e.g. \cite[Section~1]{LuZhang}); even in dimension three, this remained open until
recently \cite{Pet91}.

In a recent manuscript \cite{OAI}, OpenAI completely resolved Kobayashi's conjecture in the K\"ahler case.
The proof is closest in spirit to the quasi-negative holomorphic sectional
curvature case. For the reader's convenience, we recall the argument of \cite{WY,TY}:
\begin{enumerate}
\renewcommand{\labelenumi}{(\roman{enumi})}
\item \textcolor{red!55!black}{\emph{Volume comparison.}} Let $\w$ be a K\"ahler metric with $\Ric(\w)\ge-\w$. If $\w_0$
has holomorphic sectional curvature bounded above by $-\kappa<0$, the
Schwarz lemma applied to the identity map from $(X,\w)$ to $(X,\w_0)$ gives
\cite{WY,TY}
\[
 \w^n\ge \left(\frac{n+1}{2}\kappa\right)^n\w_0^n.
\]
\item \textcolor{red!55!black}{\emph{Standard argument.}} Since $X$ contains no rational curves, the cone theorem for
compact K\"ahler varieties \cite[Theorem~1.3]{HX} implies that $K_X$ is nef.
For a K\"ahler form $\w_0$, the theorem of Aubin and Yau \cite{Aubin,Yau} then
gives K\"ahler metrics $\w_t$ satisfying
\[
 \Ric(\w_t)=-\w_t+t\w_0,\qquad t>0.
\]
Applying the volume estimate to $\w_t$ shows that $K_X$ is big, so $X$ is Moishezon and hence projective.
The base-point-free theorem and the relative cone theorem then imply that
$K_X$ is ample, since $X$ contains no rational curves
\cite[Lemma~2.1]{DT}.
\end{enumerate}
The core of OpenAI's proof is establishing the volume comparison under hyperbolicity. Fix a Hermitian form $\w_0$ on $X$ and put
\[
 D=\sup\bigl\{|f'(0)|_{\w_0}^2:\ f\colon\D\to X\text{ holomorphic}\bigr\}.
\]
Recall Brody's criterion that
\begin{brody}[\cite{Brody}]
Let $X$ be a compact complex manifold. The following are equivalent:
\begin{enumerate}
\item $X$ is Kobayashi hyperbolic;
\item $X$ contains no entire curve, i.e., there is no nonconstant holomorphic map $\C\to X$;
\item $D<\infty$.
\end{enumerate}
\end{brody}
\noindent For every K\"ahler metric $\w$ with
$\Ric(\w)\ge-\w$, \cite[Proposition~7.1]{OAI} establishes
\begin{equation}\label{eq:oai-volume}
 \w^n\ge e^{-D}\left(\frac{n}{2n+D}\right)^n\w_0^n.
\end{equation}
To obtain this, OpenAI considers a maximizer $f_*$ of
the functional
\cite[Sections~3--5 and~7]{OAI}
\begin{equation}\label{eq:oai-functional}
 J_{q,h,\delta}(f)=|f'(0)|_{\w_0}^2
 +\frac1\pi\int_{\D}\log|z|\,f^*\w
 -\frac1{2\pi}\int_{\D}f^*\w
\end{equation}
among finite-area holomorphic discs $f$ with a fixed centre.
The second variation at $f_*$,
together with comparison with the constant disc, yields
\eqref{eq:oai-volume}.

\subsection{Volume Schwarz lemma}\label{sec:intro-volume}
We extract a general volume Schwarz lemma for extremal discs from the core argument of \cite{OAI}. 

Let's first introduce basic conventions. Let $X$ be a compact complex manifold of dimension $n$. The
Kobayashi--Royden pseudometric is defined by
\[
 F_X(x;v)=\inf\bigl\{\frac{1}{r}:\ f\colon\D\to X\text{ holomorphic},\ f(0)=x,\ f'(0)=rv\bigr\},\quad \forall x\in X,\ v\in T_{X,x}
\]
and the Kobayashi-Royden indicatrix is defined by
\[
 \mathcal I_x=\bigl\{v\in T_{X,x}:\ F_X(x;v)<1\bigr\}.
\]
For any holomorphic disc $f\colon\D\to X$ and smooth closed real $(1,1)$-form $\eta$ on $X$, the Riesz functional of $f$ with respect to $\eta$ is defined by
\begin{equation}\label{eq:R}
 R_{\eta}(f)=\frac1\pi\int_\D\log\abs z\,f^*\eta.
\end{equation}

We prove

\begin{mainthm}[Volume Schwarz lemma]\label{thm:A}
Let $(X,\w)$ be a compact K\"ahler manifold. Suppose that $\mathcal I_x$ is
bounded at some $x\in X$ and $u$ is a continuous function on $\overline{\mathcal I_x}$. Put
\begin{equation}\label{eq:J}
 J(f)=u\bigl(f'(0)\bigr)+R_\w(f).
\end{equation}
Then:
\begin{enumerate}
\renewcommand{\labelenumi}{(\alph{enumi})}
\item The functional \eqref{eq:J} attains its maximum over the holomorphic discs $f\colon\D\to X$ with $f(0)=x$, at a disc $f_*$.
\item If $u$ is of class $C^1$ near $f_*'(0)$, then $f_*$ is continuous on
$\overline\D$.
\item If $u$ is of class $C^2$ near $f_*'(0)$ and $i\partial\bar\partial u\ge0$
at $f_*'(0)$, then
\begin{equation}\label{eq:det-bound}
 \w^n(x)\ge e^{-R_{\Ric(\w)}(f_*)}\,(i\partial\bar\partial u)^n\bigl(f_*'(0)\bigr).
\end{equation}
\end{enumerate}
\end{mainthm}
\medskip
\textcolor{red!55!black}{\emph{Strategy of Theorem~\ref{thm:A}.}} (a) and (b)
will be proved in Section~\ref{sec:extremal}. For every disc $f$, by \eqref{eq:inner},
\[
 \frac1\pi\int_{\{\abs z<r\}}f^*\w\le \frac{-\Rz f}{\log\frac1r}\qquad(0<r<1),
\]
so a maximizing sequence of $J$ has bounded area on each $\{\abs z<r\}$ and no bubble at $0$, and, after passing to a subsequence, converges to a maximizer $f_*$ (Proposition~\ref{prop:exist}). A comparison argument (Proposition~\ref{prop:cap}) shows that $f_*$ is continuous on $\overline\D$.

(c) will be proved in Section \ref{sec:variations}.
Let $S$ be a bounded holomorphic matrix
function on $\D$ and $e$ be
a holomorphic frame of $f_*^*T_X$, continuous and invertible on
$\overline\D$. Deform $f_*$ with fixed centre to discs $F_\xi$ along the directions
$z\,e(z)S(z)\xi$, $\xi\in\C^n$. Lemma~\ref{lem:levi} shows
\[
 \Hess_\xi R_\w(F_\xi)\big|_{\xi=0}=-\int_{\partial\D}S^*e^*\w(f_*)eS\dd\sigma.
\]
Then $\Hess_\xi J(F_\xi)\big|_{\xi=0}\le0$ gives (Proposition~\ref{prop:variation})
\begin{equation}\label{eq:intro-variation}
 S(0)^*HS(0)\le \int_{\partial\D}S^*e^*\w(f_*)eS\dd\sigma ,\qquad
 H=e(0)^*(\Hess u)\bigl(f_*'(0)\bigr)\,e(0).
\end{equation}
After taking
determinants, the concavity of $\log\det$ shows that the best choice of $S$ should satisfy
$S^*e^*\w(f_*)eS=I$ on $\partial\D$; spectral factorization \cite{WM,HL}
provides this (see \eqref{eq:spectral}) and gives
\begin{equation}\label{eq:suffices}
 \det H\le \exp\int_{\partial\D}\log\det(e^*\w(f_*)e)\dd\sigma .
\end{equation}
Since $i\partial\bar\partial\log\det(e^*\w(f_*)e)=-f_*^*\Ric(\w)$, the
Green--Riesz formula \eqref{eq:riesz} implies (see \eqref{eq:jensen})
\begin{equation}\label{eq:intro-jensen}
 \det\bigl(e(0)^*\w(x)\,e(0)\bigr)
 =e^{-R_{\Ric(\w)}(f_*)}\exp\int_{\partial\D}\log\det(e^*\w(f_*)e)\dd\sigma ,
\end{equation}
and dividing \eqref{eq:suffices} by \eqref{eq:intro-jensen} gives
\eqref{eq:det-bound}. These steps are carried out in
Section~\ref{sec:variations}; Section~\ref{sec:proofA} proves Theorem~\ref{thm:C} and treats the equality case of Theorem~\ref{thm:B}.

\vspace{0.1cm}

\textcolor{red!55!black}{\emph{Comparison with \cite{OAI}.}}
Both \cite{OAI} and this remark compute the second variation at a maximizing disc
with fixed centre in a holomorphic frame that is unitary on the boundary.
Du--Guo--Xie used a similar frame in their earlier preprint
\cite{DGX} to construct rational curves on Fano manifolds.
\begin{enumerate}
\renewcommand{\labelenumi}{(\roman{enumi})}
\item \emph{Choice of the functional.} \cite{OAI} adds an area term to the
functional \eqref{eq:oai-functional}, because to prove that $K_X$ is nef it also maximizes
\[
 \abs{f'(0)}^2_{\w_0}-R_{\Ric(\w_0)}(f)-\frac{\delta}{2\pi}\int_\D f^*\w_0\qquad(0<\delta\le1)
\]
for a K\"ahler form $\w_0$ \cite[Proposition~6.1]{OAI}; without
the area term, the maximum might not be attained \cite[Proposition~3.2]{OAI}.
The bounded potential in \cite[Proposition~6.2]{OAI} is a special case of
Magn\'usson's Theorem~1.1 \cite{Mag12}, which \cite{OAI} does not cite.
\item \emph{Compactness.} \cite{OAI} relies global hyperbolicity while we
only control the discs $f$ with $f(0)=x$, so bubbles may form away from the
centre. Compactness modulo bubbling \cite{SU81,MS} handles this, as in
\cite{DGX}.
\item \emph{Boundary regularity.} \cite{OAI} proves Sobolev regularity of the
maximizing disc and solves a nonlinear $\bar\partial$-equation to vary it
\cite[Lemma~4.1 and Proposition~4.3]{OAI}. Once continuity up to the boundary is established
(Proposition~\ref{prop:cap}), Forstneri\v{c}'s theorem supplies the
variations \cite[Theorem~1.2]{Fo}; its Banach manifold of discs \cite[Theorem~1.1(iii)]{Fo} contains the constructions of \cite{OAI}; \cite{OAI} does not cite \cite{Fo}.
\item \emph{Second variation.} \cite{OAI} takes $u(v)=\abs v_{\w_0}^2$, and
\cite[Corollary~5.3]{OAI} is the trace of \eqref{eq:intro-variation} for
$S^*e^*\w(f_*)eS=I$ on $\partial\D$, with an extra term from the area term. We
take determinants in \eqref{eq:intro-variation} instead, which gives
\eqref{eq:suffices}.
\end{enumerate}

\subsection{Volume comparison}\label{sec:method}

If
$\Ric(\w)\ge-\w$, then $R_{\Ric(\w)}(f_*)\le-\Rz{f_*}$, and \eqref{eq:det-bound}
gives
\begin{equation}\label{eq:general-chain}
 \w^n(x)\ge e^{J(f_*)}\,\bigl(e^{-u}(i\partial\bar\partial u)^n\bigr)\bigl(f_*'(0)\bigr).
\end{equation}
Let
$\Omega\subset T_{X,x}$ be a bounded pseudoconvex domain with
$\mathcal I_x\subset\subset\Omega$.  Recall $\Omega$ carries a unique
complete K\"ahler--Einstein form $\beta_\Omega$ with
$\Ric(\beta_\Omega)=-\beta_\Omega$ \cite{ChengYau,MY}. Choose the K\"ahler-Einstein potential $u=\log(\frac{\beta_\Omega^n}{\mu}),$
where $\mu$ is a translation-invariant volume form on $T_{X,x}$. Then we have $$i\partial\bar\partial u=\beta_\Omega\
\text{ and }\ (i\partial\bar\partial u)^n=e^u\mu.$$ Since $u$ is smooth and strictly plurisubharmonic on $\Omega$, a neighbourhood of $\overline{\mathcal I_x}$, Theorem~\ref{thm:A} applies to $u$. Combine with \eqref{eq:general-chain} and
the maximality of $f_*$, we conclude that for every holomorphic disc $f\colon\D\to X$ with $f(0)=x$,
\begin{equation}\label{eq:star}
 \w^n(x)\ge e^{J(f_*)}\mu\ge e^{J(f)}\mu=e^{\Rz f}\beta_\Omega^n\bigl(f'(0)\bigr).
\end{equation}
Taking the constant disc $f$ gives \[ \w^n(x)\ge\beta_\Omega^n(0). \]
Let us introduce the intrinsic
volume form by
\[
 \mathcal S_x=
 \begin{cases}
 \sup\bigl\{\beta_\Omega^n(0):\ \Omega\subset T_{X,x}\text{ bounded
 pseudoconvex},\ \mathcal I_x\subset\subset\Omega\bigr\}, & \mathcal I_x\text{ bounded},\\[2pt]
 0, & \mathcal I_x\text{ unbounded}.
 \end{cases}
\]
Yau's volume Schwarz lemma \cite{YauSchwarz} implies $\beta_\Omega^n(0)$
decreases as $\Omega$ grows, so the larger $\mathcal I_x$, the smaller
$\mathcal S_x$. The above discussion gives

\vspace{0.1cm}

\begin{mainthm}\label{thm:B}
Let $(X,\w)$ be a compact K\"ahler manifold of dimension $n$ with
$\Ric(\w)\ge-\w$. Then
\begin{equation}\label{eq:A-domains}
 \w^n\ge \mathcal S\qquad\text{on }X .
\end{equation}
\end{mainthm}

\vspace{0.1cm}

This compares the volume of $\w$ at $x$ with the intrinsic volume form
$\mathcal S$, whereas Yau's Schwarz lemma compares it with a given K\"ahler
metric of negative holomorphic sectional curvature. Let's explain its connection with \eqref{eq:oai-volume}. One can easily check
\[
 \overline{\mathcal I_x}=\overline{\bigl\{f'(0):\ f\colon\D\to X\text{ holomorphic},\ f(0)=x\bigr\}},
\]
so
\[
 D_x(\w_0):=\sup\bigl\{\abs{f'(0)}_{\w_0}^2:\ f\colon\D\to X\text{ holomorphic},\ f(0)=x\bigr\}
 =\sup_{v\in\mathcal I_x}\abs v_{\w_0}^2.
\]
Hence, by Brody's criterion,
\[
 X\text{ is hyperbolic}\iff D_x(\w_0)<\infty\text{ for every }x\in X \iff\mathcal I_x\text{ is bounded for every }x\in X .
\]
When $D_x(\w_0)<+\infty$ at some point $x$, taking for
\[
 \Omega=\{\abs v_{\w_0}^2<R\}\ \text{with}\ R>D_x(\w_0)\quad\text{and}\quad
 \beta_\Omega^n(0)=\Bigl(\frac{n+1}{R}\Bigr)^n\w_0^n,
\]
Theorem~\ref{thm:B} gives
\begin{equation}\label{eq:Dx-volume}
 \w^n(x)\ge \Bigl(\frac{n+1}{D_x(\w_0)}\Bigr)^n\w_0^n(x),
\end{equation}
which is sharper than \eqref{eq:oai-volume} and relies only on the control at one point rather than on global hyperbolicity; and it is an equality on compact quotients of the ball with $\w_0=\w$. 

\begin{remark}
More generally, the equality in \eqref{eq:A-domains} holds on a compact quotient of some bounded symmetric domain equipped with the unique complete K\"ahler-Einstein metric.
\end{remark}

\vspace{0.1cm}

As a consequence, we obtain the following result.

\begin{mainthm}\label{thm:C}
Let $X$ be a compact K\"ahler manifold with $K_X$ nef. If the Kobayashi--Royden
indicatrix $\mathcal I_{x_0}$ is bounded at one point $x_0\in X$, then $K_X$ is
big.
\end{mainthm}

The assumption of Theorem \ref{thm:C} is strictly weaker than global hyperbolicity of $X$: it allows rational and entire curves away from
$x_0$ (see e.g. Example~\ref{ex:double-cover}), but it implies the absence of
entire curves through $x_0$. We ask
\begin{question}
  Can the nefness of
$K_X$ in Theorem \ref{thm:C} be dropped?
\end{question}
This question is related to Lang's conjecture \cite[Conjecture~5.6]{Lang86} in one direction: every subvariety of a projective hyperbolic manifold is of general type. If $X$ is hyperbolic and $\nu\colon\widetilde Y\to X$
is a resolution of singularities of a subvariety $Y$, then $\mathcal I_y$ is bounded at every point $y$
where $d\nu_y$ is injective, although $\widetilde Y$ need not be hyperbolic
or have nef canonical bundle.

\begin{remark}
   This direction of Lang's conjecture has been confirmed in the negative holomorphic sectional curvature case by Guenancia \cite{Gue}. Since $K_{\widetilde Y}$ can be expected to be pseudo-effective, Guenancia extends the volume estimate in Wu--Yau's argument to the canonical model of $(\widetilde{Y},K_{\widetilde{Y}}+tH)$ for a fixed ample divisor $H$ and $t>0$. One may check that the existence of the canonical model can be replaced by the special Fujita approximation constructed in \cite{Tru}. The natural idea is whether the holomorphic disc argument can also be carried out similarly, but it seems difficult to achieve.
\end{remark}

\vspace{0.6em}

\noindent\textbf{Acknowledgement.} The author expresses his gratitude to Westlake University and Xin Fu for supporting his research from July to November.

\vspace{0.2cm}

\noindent\textbf{AI disclosure.} The main body is written based on a talk given at a seminar after reading OpenAI's manuscript. This work was assisted by Opus 5.5 and GPT Astra-6, which were used mainly to explore the proof of Proposition \ref{prop:cap}. Other contents are formulated by the author.

\section{Preliminaries on holomorphic discs}\label{sec:prelim}
In this section we collect the basic facts on holomorphic discs used later.

\subsection*{Conventions}
We write $\D_r=\{\abs z<r\}$,
$B(a,r)=\{\abs{z-a}<r\}$, $f_r(z)=f(rz)$, and
$\mean g_r=\int_{\partial\D}g(r\zeta)\dd\sigma(\zeta)$, $\sigma$ being the normalized arc
length.

\subsection{The Riesz functional}
For $\eta=\w$, the functional $\Rz f$ in \eqref{eq:R} takes values in
$[-\infty,0]$; it is the value at the centre of the Riesz potential of the
measure $f^*\w$ \cite{Mag}. \cite{OAI} writes $P_\w=-R_\w$ and calls it the
Poisson functional. We write
\begin{equation}\label{eq:area-function}
 A_f(r,\eta)=\frac1\pi\int_{\D_r}f^*\eta,\qquad A_f(r)=A_f(r,\w).
\end{equation}

We use two elementary facts. First, since
$i\partial\bar\partial\log\abs z=\pi\delta_0$ and $\log(\frac{\abs z}{r})$ vanishes on
$\partial\D_r$, Green's formula gives the Green--Riesz representation
(see e.g. \cite[Chapter~I, (4.5)]{DemBook}, \cite[Theorem~4.5.1]{Ran})
\begin{equation}\label{eq:riesz}
 g(0)=\frac1\pi\int_{\D_r}\log\frac{\abs z}r\,i\partial\bar\partial g+\mean g_r
 \qquad\bigl(g\in C^2(\overline{\D_r})\bigr).
\end{equation}
If $\w=i\partial\bar\partial\varphi$ near $f(\overline{\D_r})$, it implies
\begin{equation}\label{eq:riesz-disc}
 \varphi(f(0))=\Rz{f_r}+\mean{\varphi\circ f}_r ,
\end{equation}
and $\Rz f=\lim_{r\to1}\Rz{f_r}$ by monotone convergence. Second, for a bounded
holomorphic map $s=\sum_{k\ge0}s_kz^k$ from $\D$ to $\C^m$,
\begin{equation}\label{eq:LP}
 \frac1\pi\int_\D\log\frac1{\abs z}\,\abs{s'}^2\,i\,dz\wedge d\bar z=\sum_{k\ge1}\abs{s_k}^2<\infty .
\end{equation}

\subsection{Compactness}
We use the compactness theory of holomorphic curves goes back to \cite{SU81} in the
following form (see e.g. \cite{MS}).

\begin{lemma}\label{lem:principal-compactness}
Let $(X,\w)$ be a compact K\"ahler manifold, and let $f_j\colon\D\to X$ be
holomorphic maps with
\[
 \sup_j A_{f_j}(r)<\infty\qquad(0<r<1).
\]
Then, after passing to a subsequence, there are a locally finite set
$\Sigma\subset\D$ and a holomorphic map $f\colon\D\to X$ with
\[
 f_j\to f\quad\text{in }C^\infty_{\mathrm{loc}}(\D\setminus\Sigma).
\]
Moreover, there is a constant $\hbar>0$ depending only on $(X,\w)$ such that
for every $a\in\Sigma$,
\[
 \liminf_{j\to\infty}\frac1\pi\int_{B(a,\rho)}f_j^*\w\ge \hbar
 \qquad\text{for every }\rho>0 .
\]
\end{lemma}
\begin{proof}
The area measures $m_j=\frac1\pi f_j^*\w$ have bounded mass on each $\D_r$, so
after passing to a subsequence they converge weakly on $\D$ to a measure $m$.
By the mean value inequality (\cite{SU81}, \cite[Lemma~4.3.1(i)]{MS}), there
are constants $c,\hbar>0$, depending only on $(X,\w)$, such that every
holomorphic map $g$ on a closed disc $\overline{B(z,r)}$ satisfies
\begin{equation}\label{eq:mean-value}
 \frac1\pi\int_{B(z,r)}g^*\w<\hbar\quad\Longrightarrow\quad
 \abs{g'(z)}_\w^2\le \frac c{r^2}\,\frac1\pi\int_{B(z,r)}g^*\w .
\end{equation}
Put
\(
 \Sigma=\bigl\{a\in\D:\ m(\{a\})\ge\hbar\bigr\}.
\)
Since $m(\D_r)<\infty$, each $\D_r$ contains only finitely many points of
$\Sigma$. If $a\notin\Sigma$, choose $\rho>0$ with $m(\overline{B(a,2\rho)})<\hbar$.
Since $\limsup_j m_j(\overline{B(a,2\rho)})\le m(\overline{B(a,2\rho)})$,
\eqref{eq:mean-value} on $B(z,\rho)\subset B(a,2\rho)$ gives, for large $j$,
\[
 \abs{f_j'(z)}_\w^2\le \frac{c\,\hbar}{\rho^2}\qquad(z\in B(a,\rho)).
\]
Hence a further subsequence converges smoothly on $\D\setminus\Sigma$. By lower
semicontinuity the limit has finite area near each point of $\Sigma$, so it
extends across $\Sigma$ by the removal of singularities
\cite[Theorem~4.1.2]{MS}. Finally, if $a\in\Sigma$ and $\rho>0$, then, since
$B(a,\rho)$ is open,
\[
 \hbar\le m(\{a\})\le m\bigl(B(a,\rho)\bigr)\le \liminf_{j\to\infty}m_j\bigl(B(a,\rho)\bigr).
\]
\end{proof}

\section{The extremal disc}\label{sec:extremal}

In this section we show that the functional \eqref{eq:J} attains its maximum
over the holomorphic discs $f\colon\D\to X$ with $f(0)=x$, and that a maximizer is continuous up to the boundary.
Throughout this section $\mathcal I_x$ is bounded, so $\overline{\mathcal I_x}$
is compact; and $u$ is a continuous function on $\overline{\mathcal I_x}$.

\subsection{Existence}\label{sec:exist}

\begin{proposition}\label{prop:exist}
The functional \eqref{eq:J} attains its maximum over the holomorphic discs $f\colon\D\to X$ with $f(0)=x$, and the
maximum is at least $u(0)$.
\end{proposition}

\begin{proof}
Since $\Rz f\le0$,
\[
 u(0)=J(\text{constant disc})\le \sup J\le \max_{\overline{\mathcal I_x}}u .
\]
Choose a maximizing sequence $(f_j)$ with $J(f_j)\ge u(0)-1$. Then
\begin{equation}\label{eq:energy-bound}
 -\Rz{f_j}\le \max_{\overline{\mathcal I_x}}u-u(0)+1 ,
\end{equation}
which is finite because $\mathcal I_x$ is bounded. For any disc $f$ and $0<r<1$, applying Fubini's theorem gives
\begin{equation}\label{eq:fubini}
 -\Rz f=\frac1\pi\int_\D\Bigl(\int_{\abs z}^1\frac{ds}s\Bigr)f^*\w
 =\int_0^1A_f(s)\,\frac{ds}s
 \ge A_f(r)\int_r^1\frac{ds}s ,
\end{equation}
where the last inequality follows from the fact that $A_f(s)$ increases with $s$; hence
\begin{equation}\label{eq:inner}
 A_f(r)\le \frac{-\Rz f}{\log\frac1r} .
\end{equation}
Combining \eqref{eq:energy-bound} and \eqref{eq:inner}, we conclude:
\begin{enumerate}
\item[(i)] for every fixed $0<r<1$ the areas $A_{f_j}(r)$ are bounded
uniformly in $j$;
\item[(ii)] $A_{f_j}(r)\to0$ as $r\to0$, uniformly in $j$.
\end{enumerate}
By (i), Lemma~\ref{lem:principal-compactness} gives a limit
$f:\D\to X$. By (ii), there is $\rho>0$ with
$\sup_jA_{f_j}(\rho)<\hbar$, so $0\notin\Sigma$ and $f_j\to f$ smoothly near
the origin,
so that $f(0)=x$ and $f_j'(0)\to f'(0)$. Since
$\abs{f_j'}^2_\w\to\abs{f'}^2_\w$ off the null set $\Sigma$, Fatou's
lemma gives $$\Rz f\ge\limsup\Rz{f_j}.$$ The continuity of $u$ gives $J(f)\ge\sup J$ and the proof is complete.
\end{proof}

\subsection{Boundary continuity}
We next show that
$f_*$ is continuous up to the boundary $\partial\D$.

\begin{proposition}\label{prop:cap}
Let $f$ maximize \eqref{eq:J} and let $u$ be of class $C^1$ near $f'(0)$. Then there is a
constant $C$ such that
\[
 \abs{f'(z)}_\w\le C\,\frac{1}{\sqrt{1-\abs z}},\qquad \forall z\in\D.
\]
In particular, $f$ extends continuously to $\overline\D$.
\end{proposition}

\begin{proof}
Fix $\zeta\in\partial\D$ and $0<\eps\le\frac14$, and let
$\Delta=\D\setminus\overline{B(\zeta,\eps)}$. Applying the Riemann mapping theorem
and Carath\'eodory's theorem gives a homeomorphism
$\phi\colon\overline\D\to\overline\Delta$, biholomorphic from $\D$ onto
$\Delta$, with $\phi(0)=0$ and $a=\phi'(0)>0$. Let
\[
 H(w)=\log\abs{\phi^{-1}(w)}-\log\abs w\qquad(w\in\overline\Delta).
\]
Then $H$ is harmonic on $\Delta$ and
continuous on $\overline\Delta$, with $H(0)=-\log a$. On $\partial\D$ we have
$H=0$. On $\partial B(\zeta,\eps)\cap\D$ we have
$\abs{\phi^{-1}(w)}=1$ and $\abs w\ge\frac34$, so
\[
 H(w)=\log\frac1{\abs w}\le\frac{1-\abs w^2}{2\abs w^2}\le 1-\abs w^2=\eps^2P_\zeta(w),
\]
where $P_\zeta(w)=\frac{1-\abs w^2}{\abs{w-\zeta}^2}$ is the Poisson kernel. The maximum principle gives $0\le H\le\eps^2P_\zeta$ on $\Delta$; at the
origin,
\begin{equation}\label{eq:cap}
 1-a \le -\log a \le \eps^2 .
\end{equation}

\emph{Comparison with $g=f\circ\phi$.} The holomorphic disc $g\colon\D\to X$ satisfies $g(0)=x$ and $g'(0)=af'(0)$. It follows from
\[
 \Rz g=\frac1\pi\int_\Delta\log\abs{\phi^{-1}(w)}\,f^*\w
\]
and $H\ge0$ that
\[
 \Rz g-\Rz f
 =\frac1\pi\int_{\D\cap B(\zeta,\eps)}\log\frac1{\abs w}\,f^*\w
 +\frac1\pi\int_\Delta H\,f^*\w
 \ge \frac1\pi\int_{\D\cap B(\zeta,\eps)}\log\frac1{\abs w}\,f^*\w .
\]
Since $J(f)\ge J(g)$, $\abs{f'(0)-af'(0)}\le\eps^2\abs{f'(0)}$ by \eqref{eq:cap}, $u$ is
$C^1$ near $f'(0)$, and $u$ is bounded on $\overline{\mathcal I_x}$, we conclude
\[
 \Rz g-\Rz f\le u\bigl(f'(0)\bigr)-u\bigl(af'(0)\bigr)\le C_1\eps^2
 \qquad(0<\eps\le\frac14).
\]
Hence,
with $C=\pi C_1$,
\begin{equation}\label{eq:weighted-cap-budget}
 \int_{\D\cap B(\zeta,\eps)}\log\frac1{\abs w}\,f^*\w\le C\eps^2
 \qquad(\zeta\in\partial\D,\ 0<\eps\le\frac14).
\end{equation}

\emph{The derivative bound.} Let $z\in\D$ with $d=1-\abs z\le\frac18$ and
$\zeta=\frac z{\abs z}$. Then $B(z,\tfrac d2)\subset\D\cap B(\zeta,2d)$, and
$\log\frac1{\abs w}\ge\frac d2$ on $B(z,\tfrac d2)$. So
\eqref{eq:weighted-cap-budget} with $\eps=2d$ gives
\[
 \frac1\pi\int_{B(z,\tfrac d2)}f^*\w\le\frac{2}{\pi d}\,C(2d)^2=\frac{8C}{\pi}\,d .
\]
If $d<\frac{\pi\hbar}{8C}$, then \eqref{eq:mean-value} with $r=\frac d2$ gives
\[
 \abs{f'(z)}_\w^2\le\frac{4c}{d^2}\cdot\frac{8C}{\pi}\,d=\frac{32cC}{\pi}\cdot\frac1{1-\abs z}.
\]
Since $\abs{f'}_\w$ is bounded on the remaining compact disc, we obtain
$\abs{f'(z)}_\w\le\frac{C_3}{\sqrt{1-\abs z}}$ on $\D$ for some constant $C_3$.

\emph{Continuity.} For $\zeta\in\partial\D$ and $0<r<s<1$, integrating along
the radius gives
\[
 \operatorname{dist}\bigl(f(r\zeta),f(s\zeta)\bigr)\le C_3\int_r^1\frac{dt}{\sqrt{1-t}}=2C_3\sqrt{1-r}.
\]
So $f(r\,\cdot)$ converges uniformly on $\partial\D$ as $r\to1$, and the proof is complete.
\end{proof}

\section{Volume Schwarz lemma}\label{sec:variations}
In this section, we complete the proof of Theorem \ref{thm:A}. Throughout, $f\colon\overline\D\to X$ is a
continuous map, holomorphic in $\D$, with $f(0)=x$ and $\Rz f>-\infty$.

\subsection{Frames along the disc}\label{sec:variation}
By applying \cite[Theorem~1.2]{Fo} to the
section $z\mapsto(z,f(z))$ of the trivial submersion $\C\times X\to\C$ over
$\overline\D$, there is a fibre-preserving biholomorphic map from a
neighbourhood of the graph of $f$ in $\C\times X$ onto an open subset of a
holomorphic vector bundle $E$ over a neighbourhood $U$ of $\overline\D$, which
maps each fibre onto a convex open set. After shrinking $U$ to a disc, $E$
is trivial (see e.g. \cite[Theorem~30.4]{Fs}). We obtain:
\begin{itemize}
\item a holomorphic map $(z,w)\mapsto G(z,w)\in X$, defined for $z\in U$ and
$w$ in a convex open set $\widetilde\Omega_z\subset\C^n$, such that
$w\mapsto G(z,w)$ is a biholomorphism onto an open subset of $X$;
\item a continuous map $a\colon\overline\D\to\C^n$, holomorphic in $\D$, with
$a(z)\in\widetilde\Omega_z$ and $G(z,a(z))=f(z)$ for $z\in\overline\D$.
\item By compactness, there is $r_0>0$ with $B(a(z),r_0)\subset\widetilde\Omega_z$
for all $z\in\overline\D$. The vertical frame
\[
 e(z)=D_wG\bigl(z,a(z)\bigr)\colon\C^n\longrightarrow T_{X,f(z)}
\]
is a holomorphic frame of $f^*T_X$ over $\D$, continuous and invertible on
$\overline\D$.
\item $e^*\w(f)e$ is continuous on $\overline\D$ and
satisfies $cI\le e^*\w(f)e\le c^{-1}I$ for some $c>0$.
\end{itemize}

\subsection{Variations and the second variation}\label{sec:second}
We state the following variation formula for an arbitrary holomorphic family.

\begin{proposition}[{cf.\ \cite[Proposition~5.2]{OAI},
\cite[Lemma~3.3]{DGX}}]\label{prop:general-variation}
Let $q$ be a smooth closed real $(1,1)$-form. Let $F\colon T\times\D_{r+\eps}\to X$
be a holomorphic map with $F_0=f$, where $T$ is a neighbourhood of $0$ in $\C^m$,
and let
\[
 V=\partial_tF\big|_{t=0}\colon\C^m\to f^*T_X
\]
be the first variation. Then
\begin{equation}\label{eq:common-hessians}
 \left.\Hess_t\left(\frac1\pi\int_{\D_r}\log\frac{|z|}{r}\,F_t^*q\right)\right|_0
 =(V^*q(f)V)(0)-\mean{V^*q(f)V}_r .
\end{equation}
\end{proposition}
\begin{proof}
Locally, we take holomorphic coordinates $w^1,\dots,w^n$ and write $q=i\partial\bar\partial\varphi$. Since $F$ is holomorphic, the chain rule gives
\[
 \partial_z\partial_{\bar z}(\varphi\circ F)
 =\sum_{a,b}\varphi_{a\bar b}(F)(\partial_zF^a)\overline{\partial_zF^b},
 \quad
 \partial_{t_k}\partial_{\bar t_j}(\varphi\circ F)
 =\sum_{a,b}\varphi_{a\bar b}(F)(\partial_{t_k}F^a)
 \overline{\partial_{t_j}F^b}.
\]
At $t=0$, the second expression is the $(j,k)$ entry of
$V^*q(f)V$. Commuting derivatives therefore gives
\begin{equation}\label{eq:commuting-variation}
 \left.\partial_{t_k}\partial_{\bar t_j}
 q(\partial_zF_t,\overline{\partial_zF_t})\right|_0
 =\partial_z\partial_{\bar z}(V^*q(f)V)_{jk}.
\end{equation}
A direct computation yields
\[
\begin{aligned}
 \left.\partial_{t_k}\partial_{\bar t_j}\right|_0
 \frac1\pi\int_{\D_r}\log\frac{|z|}{r}\,F_t^*q
 =\frac1\pi\int_{\D_r}\log\frac{|z|}{r}\,
 i\partial\bar\partial(V^*q(f)V)_{jk}=(V^*q(f)V)_{jk}(0)
 -\mean{(V^*q(f)V)_{jk}}_r,
\end{aligned}
\]
where the second equality follows from the Green--Riesz formula
\eqref{eq:riesz}.
This is the $(j,k)$ entry of the asserted identity.
\end{proof}

Let $S$ be a bounded holomorphic
$n\times n$ matrix function and, for small $\xi\in\C^n$, put
\begin{equation}\label{eq:variation}
 F_\xi(z)=G\bigl(z,a(z)+zS(z)\xi\bigr).
\end{equation}
Then $F_0=f$, $F_\xi(0)=x$, and
\begin{equation}\label{eq:jet}
 F_\xi'(0)=f'(0)+e(0)S(0)\xi,\qquad
 \partial_\xi F_\xi\big|_{\xi=0}=z\,e(z)S(z).
\end{equation}
\begin{lemma}[{cf.\ \cite[Propositions~4.3 and~5.2]{OAI}}]\label{lem:levi}
For small $\xi$, $R_\w(F_\xi)$ is finite and $C^2$. Moreover,
\begin{equation}\label{eq:riesz-boundary-hessian}
 \Hess_\xi R_\w(F_\xi)\big|_0
 =-\int_{\partial\D}S^*e^*\w(f)eS\,d\sigma.
\end{equation}
\end{lemma}
\begin{proof}
Write $F_\xi^*\w=h_\xi\,i\,dz\wedge d\bar z$. For small $\xi$, $h_\xi$ and its
$\xi$-derivatives up to order two are bounded by $C(1+|a'|+|S'|)^2$. By
\eqref{eq:LP}, this bound is integrable against $\log\frac1{|z|}$. So
$R_\w(F_\xi)$ is finite and $C^2$, and we can differentiate under the integral
sign.

Fix $0<r<1$. By \eqref{eq:jet}, $V=z\,e(z)S(z)$ and $V(0)=0$. Then
Proposition~\ref{prop:general-variation} with $q=\w$ gives
\[
 \frac1\pi\int_{\D_r}\log\frac{|z|}r\,\bigl(\Hess_\xi h_\xi\bigr)\big|_{\xi=0}\,i\,dz\wedge d\bar z
 =-r^2\int_{\partial\D}K(r\zeta)\,d\sigma(\zeta),
\]
where $K=S^*e^*\w(f)eS$. Now let $r\to1$. On $\D_r$, the integrand on the left is
bounded by $\log\frac1{|z|}\cdot C(1+|a'|+|S'|)^2$, which is integrable on $\D$.
By dominated convergence, the left side tends to
$\Hess_\xi R_\w(F_\xi)\big|_{\xi=0}$. For the right side, since $e^*\w(f)e$ is
continuous on $\overline\D$ and $S$ is bounded with radial limits almost
everywhere on $\partial\D$, we have
\[
 \abs{K(r\zeta)}\le\sup_\D\abs S^2\cdot\sup_{\overline\D}\abs{e^*\w(f)e}
 \quad\text{and}\quad
 K(r\zeta)\to K(\zeta)\ \text{ for almost every }\zeta\in\partial\D .
\]
Hence the right side tends to $-\int_{\partial\D}K\,d\sigma$, and the proof is
complete.
\end{proof}

\begin{proposition}\label{prop:variation}
Let $f$ maximize \eqref{eq:J}, $v=f'(0)$, and let $u$ be $C^2$ near
$v$. Take $e$ as in Section~\ref{sec:variation} and
$H=e(0)^*(\Hess u)(v)e(0)$. Then every bounded holomorphic matrix function $S$ on
$\D$ satisfies
\begin{equation}\label{eq:var-matrix}
 S(0)^*HS(0)\le \int_{\partial\D}S^*e^*\w(f)eS\dd\sigma .
\end{equation}
\end{proposition}

\begin{proof}
By Proposition~\ref{prop:cap}, $f$ is continuous on $\overline\D$, so
Section~\ref{sec:variation} applies. By \eqref{eq:jet} and
Lemma~\ref{lem:levi}, the function
\[
 \xi\longmapsto J(F_\xi)=u\bigl(v+e(0)S(0)\xi\bigr)+\Rz{F_\xi}
\]
is $C^2$ near $0$. Since $F_\xi(0)=x$ and $f$ maximizes $J$, it attains its
maximum at $\xi=0$, hence
\[
 \Hess_\xi J(F_\xi)\big|_{\xi=0}
 =S(0)^*HS(0)-\int_{\partial\D}S^*e^*\w(f)eS\dd\sigma\le0 .
\qedhere
\]
\end{proof}

\subsection{Proof of Theorem~\ref{thm:A}}

Part (a) is Proposition~\ref{prop:exist} and part (b) is
Proposition~\ref{prop:cap}. For (c), put $f=f_*$ and $v=f_*'(0)$, and assume
$(\Hess u)(v)\ge0$. In the frame $e$ of Section~\ref{sec:variation},
\[
 i\partial\bar\partial\log\det(e^*\w(f)e)=-f^*\Ric(\w).
\]
Since $\log\det(e^*\w(f)e)$ is continuous on $\overline\D$ and
$R_{\Ric(\w)}(f)$ is finite, the Green--Riesz formula \eqref{eq:riesz} on
$\D_r$ gives, as $r\to1$,
\begin{equation}\label{eq:jensen}
 \log\det(e(0)^*\w(x)e(0))=\int_{\partial\D}\log\det(e^*\w(f)e)\,d\sigma
 -R_{\Ric(\w)}(f).
\end{equation}

Let $S$ be a bounded holomorphic matrix function with bounded inverse, and put
$K=S^*e^*\w(f)eS$. Taking determinants in \eqref{eq:var-matrix},
\[
 \frac{\det(\Hess u)(v)}{\det\w(x)}
 =\frac{\det(S(0)^*HS(0))}{\det(S(0)^*e(0)^*\w(x)e(0)S(0))}
 \le\frac{\det\int_{\partial\D}K\,d\sigma}{\det(S(0)^*e(0)^*\w(x)e(0)S(0))}.
\]
Since $\log\abs{\det S}^2$ is bounded and harmonic, \eqref{eq:jensen} gives
\[
 \log\det(S(0)^*e(0)^*\w(x)e(0)S(0))
 =\int_{\partial\D}\log\det K\,d\sigma-R_{\Ric(\w)}(f),
\]
and therefore
\begin{equation}\label{eq:determinant-defect}
 \frac{\det(\Hess u)(v)}{\det\w(x)}
 \le e^{R_{\Ric(\w)}(f)}
 \exp\Bigl(\log\det\int_{\partial\D}K\,d\sigma-\int_{\partial\D}\log\det K\,d\sigma\Bigr).
\end{equation}
By the concavity of $\log\det$, the last factor is at least $1$, with equality
when $K=I$ on $\partial\D$.

It remains to find $S$ with $K=I$ on $\partial\D$. Since $e^*\w(f)e$ and its
inverse are bounded, the theorem of Wiener--Masani and Helson--Lowdenslager
\cite[Theorem~7.13]{WM}, \cite{HL} gives an outer matrix function $\Phi$ with
\begin{equation}\label{eq:spectral}
 e^*\w(f)e=\Phi^*\Phi\quad\text{almost everywhere on }\partial\D .
\end{equation}
The upper bound on $e^*\w(f)e$ gives $\Phi\in H^\infty$. Since $\det\Phi$ is
outer and $\abs{\det\Phi}\ge c^{\frac n2}$ almost everywhere on $\partial\D$,
\(
 \Phi^{-1}=\frac{\operatorname{adj}\Phi}{\det\Phi}\in H^\infty .
\)
Taking $S=\Phi^{-1}$ gives $K=I$ on $\partial\D$, and \eqref{eq:determinant-defect}
becomes \eqref{eq:det-bound}.

\section{The volume inequalities}\label{sec:proofA}
\subsection{Proof of Theorem~\ref{thm:B}}
Theorem~\ref{thm:B} follows from Theorem~\ref{thm:A} as in Section~\ref{sec:method}.

The inequality \eqref{eq:star} is sharp on compact quotients $X=\Gamma\backslash D$
of bounded symmetric domains, with $\w$ induced by $\beta_D$. Realize $D$ in
Harish-Chandra coordinates, with $x$ the image of $0$, and consider the linear
discs $f(\zeta)=\zeta v$, $v\in D$. The potential $\psi=\log\det g_D$ of
$\beta_D$, where $g_D$ is the coefficient matrix of $\beta_D$ in these coordinates, is invariant under the rotations $v\mapsto e^{i\theta}v$, which are
linear automorphisms of $D$ with Jacobian of modulus one, so
\eqref{eq:riesz-disc} gives $\Rz f=\psi(0)-\psi(v)$, that is,
\[
 e^{\Rz f}\beta_D^n(v)=\beta_D^n(0)=\w^n(x).
\]
By \cite[Example~1.17]{JP91}, the Kobayashi--Royden metric of $D$ at $0$ is its
Minkowski function, as $D$ is bounded, convex and balanced. Every holomorphic disc $h\colon\D\to X$ with $h(0)=x$ lifts to a holomorphic disc $\tilde h\colon\D\to D$ with $\tilde h(0)=0$. Hence
\[
 \mathcal I_x=D .
\]
Since $\beta_\Omega^n(0)$ is decreasing in $\Omega$ \cite{YauSchwarz} and
$\beta_{rD}^n(0)=r^{-2n}\beta_D^n(0)$, we get
\[
 \mathcal S_x=\beta_D^n(0)=\w^n(x).
\]
So \eqref{eq:A-domains} is an equality. Now let $g\colon\D\to X$ be a holomorphic disc with $g(0)=x$ and $g'(0)=v$. Then \eqref{eq:star} with $\Omega=rD$, $r>1$, gives
$e^{\Rz g}\beta_{rD}^n(v)\le\w^n(x)$, and since
$\beta_{rD}^n(v)=r^{-2n}\beta_D^n(\tfrac vr)\to\beta_D^n(v)$ as $r\downarrow1$,
\[
 \Rz g\le\Rz f .
\]
Therefore the linear disc maximizes $R_\w$ among the holomorphic discs $g\colon\D\to X$ with $g(0)=x$ and $g'(0)=v$, and equality holds in \eqref{eq:star} for it.

\subsection{Bigness}\label{sec:bigness}
Now we use Theorem~\ref{thm:B} to prove Theorem~\ref{thm:C}. The argument is
standard, and we include it for the reader's convenience.

\begin{proof}[Proof of Theorem~\ref{thm:C}]
Let $\w_0$ be a K\"ahler form. Since $K_X$ is nef, the class
$2\pi c_1(K_X)+t[\w_0]$ is K\"ahler for every $t>0$. By the theorem of Aubin and
Yau \cite{Aubin,Yau}, it contains a K\"ahler form $\w_t$ satisfying
\[
 \Ric(\w_t)=-\w_t+t\,\w_0 .
\]
So $-\w_t\le\Ric(\w_t)\le\w_0$ when $0<t\le1$.

Take a bounded pseudoconvex domain
$\Omega\subset T_{X,x_0}$ with $\mathcal I_{x_0}\subset\subset\Omega$. Then
Theorem~\ref{thm:B} gives
\[
 \w_t^n(x_0)\ge\beta_\Omega^n(0)>0\qquad(0<t\le1).
\]
Take a coordinate ball $B$
centred at $x_0$ with $\w_0=i\partial\bar\partial\rho$ on $\overline B$, and
write $\w_t^n=q_t\,dV$, where $dV$ is the Euclidean volume form. Then
\[
 -i\partial\bar\partial\log q_t=\Ric(\w_t)\le i\partial\bar\partial\rho ,
\]
so $\log q_t+\rho$ is plurisubharmonic on $B$, and so is $q_te^\rho$. The sub-mean
value inequality at $x_0$ gives
\[
 \int_X\w_t^n\ge\int_Bq_t\dd V\ge\Vol(B)\,e^{\rho(x_0)-\sup_B\rho}\,q_t(x_0)\ge c_2
\]
for a constant $c_2>0$ independent of $t$. Letting $t\to0$, we get
\[
 (2\pi)^n\int_Xc_1(K_X)^n\ge c_2>0 .
\]
So $K_X$ is big and the proof is complete.
\end{proof}

\begin{example}\label{ex:double-cover}
Let $M$ be a Kobayashi hyperbolic projective
manifold, $Y\subset M$ an irreducible normal surface, and
$\pi\colon X\to Y$ its minimal resolution. Every rational curve on $X$ is
contracted by $\pi$, since $M$ contains no rational curves. Thus $X$ is not
ruled, and the minimality of $\pi$ implies that $X$ contains no $(-1)$-curves.
By the classification of surfaces, $K_X$ is nef.

At any $x\in\pi^{-1}(Y_{\mathrm{reg}})$ the differential of $X\to M$ is
injective. The distance-decreasing property of the Kobayashi--Royden
pseudometric therefore implies $\mathcal I_x$ bounded. Theorem~\ref{thm:C}
gives that $K_X$ is big, so $Y$ is of general type. If $Y$ has an ordinary
double point, the corresponding exceptional curve is a rational $(-2)$-curve;
thus $X$ need not be hyperbolic.
\end{example}


\begin{thebibliography}{WY16b}\small
\bibitem[Aub78]{Aubin} T.~Aubin, \emph{\'Equations du type Monge--Amp\`ere sur les
vari\'et\'es k\"ahl\'eriennes compactes}, Bull. Sci. Math. \textbf{102} (1978),
63--95.
\bibitem[Bro78]{Brody} R.~Brody, \emph{Compact manifolds and hyperbolicity}, Trans.
Amer. Math. Soc. \textbf{235} (1978), 213--219.
\bibitem[Cam91]{Cam91} F.~Campana, \emph{Twistor spaces and nonhyperbolicity of
certain symplectic K\"ahler manifolds}, in: Complex Analysis (Wuppertal,
1991), K.~Diederich (ed.), Aspects Math. \textbf{E17}, Vieweg, 1991, 64--69.
\bibitem[CY18]{CY} B.-L.~Chen and X.~Yang, \emph{Compact K\"ahler manifolds
homotopic to negatively curved Riemannian manifolds}, Math. Ann. \textbf{370}
(2018), 1477--1489.
\bibitem[CY80]{ChengYau} S.-Y.~Cheng and S.-T.~Yau, \emph{On the existence of a complete
K\"ahler metric on non-compact complex manifolds and the regularity of
Fefferman's equation}, Comm. Pure Appl. Math. \textbf{33} (1980), 507--544.
\bibitem[Dem]{DemBook} J.-P.~Demailly, \emph{Complex Analytic and Differential
Geometry}, book available at https://www-fourier.ujf-grenoble.fr/\~{}demailly/books.html.
\bibitem[Dem85]{Dem} J.-P.~Demailly, \emph{Champs magn\'etiques et in\'egalit\'es de
Morse pour la $d''$-cohomologie}, Ann. Inst. Fourier \textbf{35} (1985),
189--229.
\bibitem[DP04]{DP} J.-P.~Demailly and M.~P\u{a}un, \emph{Numerical characterization of
the K\"ahler cone of a compact K\"ahler manifold}, Ann. of Math. \textbf{159}
(2004), 1247--1274.
\bibitem[Div20]{DiverioSurvey} S.~Diverio, \emph{Kobayashi hyperbolicity,
negativity of the curvature and positivity of the canonical bundle},
arXiv:2011.11379v2 (2020).
\bibitem[DT19]{DT} S.~Diverio and S.~Trapani, \emph{Quasi-negative holomorphic
sectional curvature and positivity of the canonical bundle}, J. Differential
Geom. \textbf{111} (2019), 303--314.
\bibitem[DGX26]{DGX} Y.-H.~Du, B.~Guo and S.-Y.~Xie, \emph{Analytic construction of
rational curves on Fano manifolds}, arXiv:2609.11612v2 (30 September 2026).
\bibitem[For81]{Fs} O.~Forster, \emph{Lectures on Riemann Surfaces}, Grad. Texts in
Math. \textbf{81}, Springer, 1981.
\bibitem[For07]{Fo} F.~Forstneri\v{c}, \emph{Manifolds of holomorphic mappings from
strongly pseudoconvex domains}, Asian J. Math. \textbf{11} (2007), 113--126.
\bibitem[Gue22]{Gue} H.~Guenancia, \emph{Quasiprojective manifolds with negative
holomorphic sectional curvature}, Duke Math. J. \textbf{171} (2022), no.~2,
417--442.
\bibitem[HX26]{HX} C.~Hacon and L.~Xie, \emph{On the K\"ahler MMP and the
transcendental base-point-free theorem}, preprint, arXiv:2607.24986v1 (2026).
\bibitem[HL61]{HL} H.~Helson and D.~Lowdenslager, \emph{Prediction theory and Fourier
series in several variables. II}, Acta Math. \textbf{106} (1961), 175--213.
\bibitem[JP91]{JP91} M.~Jarnicki and P.~Pflug, \emph{Invariant pseudodistances and
pseudometrics---completeness and product property}, Ann. Polon. Math.
\textbf{55} (1991), 169--189.
\bibitem[Kob70]{Kob} S.~Kobayashi, \emph{Hyperbolic Manifolds and Holomorphic
Mappings}, Pure Appl. Math. \textbf{2}, Marcel Dekker, New York, 1970;
2nd ed., World Scientific, 2005.
\bibitem[Kob76]{Kob76} S.~Kobayashi, \emph{Intrinsic distances, measures and geometric
function theory}, Bull. Amer. Math. Soc. \textbf{82} (1976), 357--416.
\bibitem[Lan86]{Lang86} S.~Lang, \emph{Hyperbolic and Diophantine analysis}, Bull.
Amer. Math. Soc. (N.S.) \textbf{14} (1986), 159--205.
\bibitem[LZ17]{LuZhang} S.~S.~Y.~Lu and D.-Q.~Zhang, \emph{Positivity criteria for
log canonical divisors and hyperbolicity}, J. Reine Angew. Math. \textbf{726}
(2017), 173--186.
\bibitem[Mag11]{Mag} B.~S.~Magn\'usson, \emph{Extremal $\w$-plurisubharmonic
functions as envelopes of disc functionals}, Ark. Mat. \textbf{49} (2011),
383--399.
\bibitem[Mag12]{Mag12} B.~S.~Magn\'usson, \emph{Extremal $\w$-plurisubharmonic
functions as envelopes of disc functionals: generalization and applications
to the local theory}, Math. Scand. \textbf{111} (2012), 296--319.
\bibitem[MS12]{MS} D.~McDuff and D.~Salamon, \emph{$J$-holomorphic Curves and
Symplectic Topology}, 2nd ed., Amer. Math. Soc. Colloq. Publ. \textbf{52}, 2012.
\bibitem[MY83]{MY} N.~Mok and S.-T.~Yau, \emph{Completeness of the K\"ahler--Einstein
metric on bounded domains and the characterization of domains of holomorphy
by curvature conditions}, in: The Mathematical Heritage of Henri Poincar\'e,
Part~1, Proc. Sympos. Pure Math. \textbf{39}, Amer. Math. Soc., 1983, 41--60.
\bibitem[OAI26]{OAI} OpenAI, \emph{Canonical ampleness of compact hyperbolic K\"ahler
manifolds}, preprint, 23 September 2026, \url{https://github.com/openai/math}
(directory \texttt{preprints}).
\bibitem[Pet91]{Pet91} T.~Peternell, \emph{Calabi--Yau manifolds and a conjecture of
Kobayashi}, Math. Z. \textbf{207} (1991), 305--318.
\bibitem[Ran95]{Ran} T.~Ransford, \emph{Potential Theory in the Complex Plane}, London
Math. Soc. Stud. Texts \textbf{28}, Cambridge Univ. Press, 1995.
\bibitem[SU81]{SU81} J.~Sacks and K.~Uhlenbeck, \emph{The existence of minimal
immersions of $2$-spheres}, Ann. of Math. (2) \textbf{113} (1981), 1--24.
\bibitem[TY17]{TY} V.~Tosatti and X.~Yang, \emph{An extension of a theorem of Wu--Yau},
J. Differential Geom. \textbf{107} (2017), 573--579.
\bibitem[Tru26]{Tru} A.~Trusiani, with an appendix by S.~Boucksom, \emph{A solution to the
Yau--Tian--Donaldson conjecture through special Fujita approximations},
arXiv:2605.30063 (2026).
\bibitem[WM57]{WM} N.~Wiener and P.~Masani, \emph{The prediction theory of
multivariate stochastic processes. I}, Acta Math. \textbf{98} (1957), 111--150.
\bibitem[WY16a]{WY16remark} D.~Wu and S.-T.~Yau, \emph{A remark on our paper
``Negative holomorphic curvature and positive canonical bundle''}, Comm. Anal.
Geom. \textbf{24} (2016), 901--912.
\bibitem[WY16b]{WY} D.~Wu and S.-T.~Yau, \emph{Negative holomorphic curvature and
positive canonical bundle}, Invent. Math. \textbf{204} (2016), 595--604.
\bibitem[Yau78a]{YauSchwarz} S.-T.~Yau, \emph{A general Schwarz lemma for K\"ahler
manifolds}, Amer. J. Math. \textbf{100} (1978), 197--203.
\bibitem[Yau78b]{Yau} S.-T.~Yau, \emph{On the Ricci curvature of a compact K\"ahler
manifold and the complex Monge--Amp\`ere equation. I}, Comm. Pure Appl. Math.
\textbf{31} (1978), 339--411.
\end{thebibliography}
\end{document}